\documentclass{article}

\usepackage{amsmath}
\usepackage{amssymb}
\usepackage{amsfonts}
\usepackage{latexsym}
\usepackage{enumerate}

\newtheorem{llemma}{Lemma}[section]
\newtheorem{prop}[llemma]{Proposition}

\newtheorem{exmp}[llemma]{Example}
\newtheorem{ttheorem}[llemma]{Theorem}
\newtheorem{ccorollary}[llemma]{Corollary}
\newtheorem{defn}[llemma]{Definition}

\newtheorem{key}[llemma]{Keyword}

\begin{document}

\title{{ {Application of Soft Sets in Decision Making based on Game Theory}}}
\author{Irfan Deli, Naim \c{C}a\u{g}man \\ \\
         Department of Mathematics, Faculty of Arts and
          Sciences,\\
 Kilis 7 Aral{\i}k University, 79000 Kilis, Turkey, \\
              irfandeli@kilis.edu.tr \\~~\\
             Department of Mathematics, Faculty of Arts and
          Sciences,\\
Gaziosmanpa\c{s}a University, 60250 Tokat, Turkey \\
             naim.cagman@gop.edu.tr
       }
\date{}

\maketitle


\begin{abstract}
In this work, after given the definition of soft sets and their
basic operations we define two person soft games which can apply to
problems contain vagueness and uncertainty. We then give four
solution methods of the games which are soft saddle points, soft
lover and soft upper values, soft dominated strategies and soft Nash
equilibrium. We also give an example from the real world which shows
that the methods can be successfully applied to a financial problem.
Finally, we extended the two person soft games to n-person soft
games.
\begin{key}
Soft sets, two person soft games, soft payoff functions, soft
dominated strategies, soft lover and soft upper values, soft Nash
equilibrium.
\end{key}
\end{abstract}


\section{Introduction}

In 1999, Molodtsov \cite{mol-99} introduced soft set theory for
modeling vagueness and uncertainty. In \cite{mol-99}, Molodtsov
pointed out several directions for the applications of soft sets,
such as stability and regularization, game theory and operations
research, and soft analysis. After Molodtsov, works on soft set
theory has been progressing rapidly. For instance; on the theory of
soft sets \cite{ali-09,cag-09a,cag-09b,maj-03,mol-01,mol-04,mol-06},
on the soft decision making \cite{deli-12a,deli-12b,deli-12,maj-02},
on the algebraic structures soft sets  \cite{akt-07,op1,op2,op3} are
some of the selected works.

Game theory is originally the mathematical study of competition and
cooperation. In other words, game theory is a study of strategic
decision making \cite{mye-91}. Game theory was introduced in 1944
with the publication of von Neumann and Morgenstern \cite{neo-44}.
They started modern game theory with the two-person zero-sum games
and its proof. Game theory is mainly used in many fields such as;
economics, political science, psychology and so on \cite{aum-x}.
Ferguson \cite{fer-08} present various mathematical models of game
theory. Binmore \cite{bin-82} focused the cooperative and
noncooperative game theory. Aliprantis and Chakrabarti
\cite{alip-00} give games with decision making.

In 1965, Zadeh \cite{zad-65} developed the theory of fuzzy sets that
is the most appropriate theory for dealing with uncertainties. In
recent years, many interesting applications of game theory have been
expanded by embedding the ideas of fuzzy sets. The two person
zero-sum games with fuzzy payoffs and fuzzy goals game theory have
been studied by many authors (e.g. \cite{bec-04, cev-09,cev-10,
che-06, lar-09, mae-03, mae-00, vij-05}). The max-min solution with
respect to a degree of attainment of a fuzzy goal has also been
studied (e.g. \cite{abb-06, dex-04,sak-94,rus-00,yeh-08,yos-99}).
Many study of game theory have been expanded by using the ideas of
interval value (e.g. \cite{dex-04, liu-09, mal-11}). Theory for
linear programming problems with fuzzy parameters is introduced
(e.g. \cite{bec-04, bec-05, cam-89}).

In the classical and fuzzy games, the payoff functions are reel
valued, and therefore the solution of such games are obtained by
using arithmetic operations. Especially, fuzzy games depend on the
fuzzy set that is described by its membership function. It is
mentioned in \cite{mol-99}, there exists a difficulty to set the
membership function in each particular case, and also the fuzzy set
operations based on the arithmetic operations with membership
functions do not look natural since the nature of the membership
function is extremely individual.

In this work, we propose a game model for dealing with uncertainties
which is free of the difficulties mentioned above. The proposed new
game is called a soft game since it is based on soft sets theory. To
construct a soft set we can use any parametrization with the help of
words and sentences, real numbers, functions, mappings, and so on.
Therefore, payoff functions of the soft game are set valued function
and solution of the soft games obtained by using the operations of
sets that make this game very convenient and easily applicable in
practice. The present expository paper is a condensation of part of
the dissertation \cite{deli-13tez}.

This work is organized as follows. In the next section, most of the
fundamental definitions of the operations of soft sets are
presented. In Section 3, we construct two person soft games and then
give four solution methods for the games which are soft saddle
points, soft lover and soft upper value, soft dominated strategies
and soft Nash equilibrium. In section 4, we give an application for
two person soft games. In section 5, we give n-person soft games
that is extension of the two person  soft games. In final Section,
we concluded the work.
\section{Soft sets}\label{ss}

In this section, we present the basic definitions and results of
soft set theory \cite{cag-09a}. More detailed explanations related
to this subsection may be found in \cite{cag-09a, maj-03, mol-99}.

\begin{defn}\cite{mol-99}
Let $U$ be a universe, $P(U)$ be the power set of $U$ and $E$ be a
set of parameters that are describe the elements of $U$. A soft set
$S$ over $U$ is a set defined by a set valued function $S$
representing a mapping
$$
f_S: E\to P(U)
$$
It is noting that the soft set is a parametrized family of subsets
of the set $U$, and therefore it can be written a set of ordered
pairs
$$
S= \{(x, f_S(x)): x\in E\}
$$
\end{defn}
Here, $f_S$ is called approximate function of the soft set $S$ and
$f_S(x)$ is called $x$-approximate value of $x \in E$. The subscript
$S$ in the $f_S$ indicates that $f_S$ is the approximate function of
$S$.

Generally, $f_S$, $f_T$, $f_V$, ... will be used as an approximate
functions of $S$, $T$, $V$, ..., respectively.

Note that if $f_S(x)=\emptyset$, then the element $(x, f_S(x))$ is
not appeared in $S$.
\begin{exmp}\label{ex-soft}
Suppose that $U=\{u_1,u_2,u_3,u_4\}$ is the universe contains four
cars under consideration in an auto agent and
$E=\{x_1,x_2,x_3,x_4\}$ is the set of parameters, where $x_i$
$(i=1,2,3,4)$ stand for `safety'', ``cheap'', ``modern'' and
``large'', respectively.

A customer to select a car from the auto agent, can construct a soft
set $S$ that describes the characteristic of cars according to own
requests. Assume that $f_S(x_1)= \{u_1, u_2\}$, $f_S(x_2)= \{u_1,
u_2,u_4 \}$, $ f_S(x_3)= \emptyset$, $ f_S(x_4)= U$ then the
soft-set $S$ is written by
$$
S=\{(x_1, \{u_1, u_2\}),(x_2, \{u_1, u_2,u_4\}), (x_4,U)\}
$$

By using same parameter set $E$, another customer to select a car
from the same auto agent, can construct a soft set $T$ according to
own requests. Here $T$ may be different then $S$. Assume that
$f_T(x_1)= \{u_1, u_2\}$, $f_T(x_2)= \{u_1, u_2,u_3 \}$, $
f_T(x_3)=\{u_1, u_2\}$, $ f_T(x_4)= \{u_1\}$ then the soft-set $T$
is written by
$$
T=\{(x_1, \{u_1, u_2\}), (x_2, \{u_1, u_2,u_3\}), (x_3, \{u_1,
u_2\}), (x_4,\{u_1\}) \}
$$

\end{exmp}

\begin{defn}\cite{cag-09a} Let $S$ and $T$ be two soft sets. Then,
\begin{enumerate}[a)]
    \item
If $f_S(x)=\emptyset$ for all $x\in E$, then $S$ is called a empty
soft set, denoted by $S_\Phi$.
    \item
If $f_S(x)\subseteq f_T(x)$ for all $x\in E$, then $S$ is a soft
subset of $T$, denoted by $S \tilde{\subseteq} T$.
\end{enumerate}
\end{defn}

\begin{defn}\cite{cag-09a} Let $S$ and $T$ be two soft sets. Then,
\begin{enumerate}[a)]
\item
Complement of $S$ is denoted by $S^{\tilde c}$. Its approximate
function $f_{S^{\tilde c}}$ is defined by
$$
f_{S^{\tilde c}}(x)= U\setminus f_S(x)\quad\textrm{for all } x \in E
$$
\item
Union of $ S$ and $T$ is denoted by $S \tilde{\cup} T$. Its
approximate function $f_{S \widetilde{\cup} T}$ is defined by
$$
f_{S \widetilde{\cup} T}(x)= f_S(x)\cup f_T(x)\quad\textrm{for all }
x \in E.
$$

\item Intersection of $S$ and $T$ is
denoted by $S\tilde{\cap}T$. Its approximate function $f_{S
\widetilde{\cap} T}$ is defined by
$$
f_{S \widetilde{\cap} T}(x)= f_S(x)\cap f_T(x)\quad\textrm{for all }
x \in E.
$$
\end{enumerate}
\end{defn}
\section{Two Person Soft Games}\label{ss}
In this section,  we construct two person  soft games with soft
payoffs. We then give four solution methods for the games. The basic
definitions and preliminaries of the game theory we refer to
\cite{alip-00,fer-08,mol-99,neo-44,owe-95,webb-88}.

\begin{defn} Let $X,Y$ are a sets of strategies.
A choice of behavior is called an action. The elements of $X\times
Y$ are called action pairs. That is, $X\times Y$ is the set of
available actions.
\end{defn}
\begin{defn}
Let $U$ be a set of alternatives,  $P(U)$ be the power set of $U$,
$X,Y$ are a sets of strategies. Then, a set valued function
$$
f_S: X\times Y \rightarrow P(U)
$$
is called a soft payoff function. For each $(x,y) \in X\times Y$,
the value $f_S(x,y)$ is called a soft payoff.
\end{defn}

\begin{defn}
Let $X$ and $Y$ be a set of strategies of Player 1 and 2,
respectively, $U$ be  a set of alternatives and $f_{S_k}:X\times
Y\rightarrow P(U)$ be a soft payoff function for player $k$,
$(k=1,2)$. Then, for each Player $k$, a two person soft game
($tps$-game) is defined by a soft set over $U$ as
$$
{S_k}=\{((x,y),f_{S_k}(x,y)): (x,y)\in X\times Y\}
$$
The $tps$-game is played as follows: at a certain time Player 1
chooses a strategy $x_i\in X$, simultaneously Player 2 chooses a
strategy $y_j\in Y$ and once this is done each player k (k=1,2)
receives the soft payoff $f_{S_k}(x_i,y_j)$.

If $X=\{x_1,x_2,...,x_m\}$ and $Y=\{y_1,y_2,...,y_n\}$, then the
soft payoffs of ${S_k}$ can be arranged in the form of the $m\times
n$ matrix shown in Table 1.
$$
\begin{tabular}{|c|c|c|c|c|}
\hline
\({S_k}\)& \(y_1\)                  & \(y_2\)                  &...&\(y_n\) \\
\hline
\(x_1\)  & \(f_{S_k}(x_1,y_1)\)&  \(f_{S_k}(x_1,y_2)\) &...& \( f_{S_k}(x_1,y_n)\) \\
\hline
\(x_2\)  & \(f_{S_k}(x_2,y_1)\)&  \(f_{S_k}(x_2,y_2)\) &...& \( f_{S_k}(x_2,y_n)\) \\
\hline
\vdots   & \vdots                     &  \vdots                &\(\ddots\)&  \vdots \\
\hline
\(x_m\)  & \(f_{S_k}(x_m,y_1)\)&  \(f_{S_k}(x_m,y_2)\) &...&  \(f_{S_k}(x_m,y_n)\) \\
\hline
\end{tabular}
$$
\begin{center}
\footnotesize{\emph{\emph{Table 1}: The two person soft game}}
\end{center}
\end{defn}

Now, we can give an example for $tps$-game.

\begin{exmp}\label{1}Let $U=
\{u_1,u_2,u_3, u_4,u_5,u_6, u_7, u_8,u_9,u_{10}\}$be  a set of
alternatives,  $P(U)$ be the power set of $U$, $X=\{x_1, x_3, x_5\}$
and $Y=\{x_1, x_2, x_4\}$ be a set of the strategies Player 1 and 2,
respectively.

If Player 1 constructs a $tps$-games as follows,
$$
\begin{array}{rr}
{S_1}=& \bigg\{ ((x_1,x_1),\{u_1,u_2,u_5,u_8\})
,(x_1,x_2),\{u_1,u_2,u_3,u_4,u_5,u_8\}),(x_1,x_4),\\&\{u_3,u_8\}),((x_3,x_1),\{u_1,u_3,u_7\})
,(x_3,x_2),\{u_1,u_2,u_3,u_5,u_6,u_7\}),\\&(x_3,x_4),\{u_1,u_2,u_3\}),((x_5,x_1),\{u_3,u_4,u_5,u_8\}
)
,(x_5,x_2),\{u_1,u_2,u_3,\\&u_4,u_5,u_6,u_8\}),(x_5,x_4),\{u_1,u_2,u_3,u_8\})\bigg\} \\
\end{array}
$$
then the soft payoffs of the game can be arranged as in Table 2,
$$
\begin{tabular}{|c|c|c|c|c|c|c|c|}
\hline
${S_1}$& \(x_{1}\) & \(x_{2}\) &  \(x_{4}\)\\
\hline
\(x_1\) & \(\{u_1,u_2,u_5,u_8\} \) & \(\{u_1,u_2,u_3,u_4,u_5,u_8\} \) &  \(\{u_3,u_8\} \) \\
\hline
\(x_3\) & \(\{u_1,u_3,u_7\} \) & \(\{u_1,u_2,u_3,u_5,u_6,u_7\} \) &  \(\{u_1,u_2,u_3\}\) \\
\hline
\(x_5\) & \(\{u_3,u_4,u_5,u_8\} \) & \(\{u_1,u_2,u_3,u_4,u_5,u_6,u_8\} \) &  \(\{u_1,u_2,u_3,u_8\}\) \\
\hline
\end{tabular}
$$
\begin{center}
\footnotesize{\emph{Table 2}}
\end{center}

Let us explain some element of this game; if Player 1 select $x_3$
and Player 2 select $x_2$, then the value of game will be a set
$\{u_1,u_2,u_3,u_5,u_6,u_7\}$, that is,
$$
f_{S_1}(x_3,x_2)=\{u_1,u_2,u_3,u_5,u_6,u_7\}
$$
In this case, Player 1 wins the set of alternatives
$\{u_1,u_2,u_3,u_5,u_6,u_7\}$ and Player 2 lost the same set of
alternatives.

Similarly, if Player 2 constructs a $tps$-game as follows,
$$
\begin{array}{rl}
{S_2}=& \bigg\{ ((x_1,x_1),\{u_3,u_4,u_6,u_7\} )
,(x_1,x_2),\{u_6,u_7\} ),(x_1,x_4),\{u_1,u_2,u_4,\\&u_5,u_6,u_7\}
),((x_3,x_1),\{u_2,u_4,u_5,u_6,u_8\}) ,(x_3,x_2),\{u_4,u_8\}
),(x_3,x_4),\\&\{u_4,u_5,u_6,u_7,u_8\}),((x_5,x_1),\{u_1,u_2,u_6u_7\}
)
,(x_5,x_2),\{u_7\} ),\\&(x_5,x_4),\{u_4,u_5,u_6,u_7\})\bigg\} \\
\end{array}
$$
then the soft payoffs of the game can be arranged as in Table 3,

$$
\begin{tabular}{|c|c|c|c|c|c|c|c|}
\hline
${S_2}$& \(x_{1}\) & \(x_{2}\) &  \(x_{4}\)\\
\hline
\(x_1\) & \(\{u_3,u_4,u_6,u_7\} \) & \(\{u_6,u_7\} \) &  \(\{u_1,u_2,u_4,u_5,u_6,u_7\} \) \\
\hline
\(x_3\) & \(\{u_2,u_4,u_5,u_6,u_8\} \) & \(\{u_4,u_8\} \) &  \(\{u_4,u_5,u_6,u_7,u_8\}\) \\
\hline
\(x_5\) & \(\{u_1,u_2,u_6,u_7\} \) & \(\{u_7\} \) &  \(\{u_4,u_5,u_6,u_7\}\) \\
\hline
\end{tabular}
$$
\begin{center}
\footnotesize{\emph{Table 3}}
\end{center}

Let us explain some element of this $tps$-game; if Player 1 select
$x_3$ and Player 2 select $x_2$, then the value of game will be a
set $\{u_4,u_8\}$, that is,
$$
f_{S_2}(x_3,x_2)=\{u_4,u_8\}
$$

In this case, Player 1 wins the set of alternatives $\{u_4,u_8\}$
and Player 2 lost $\{u_4,u_8\}$.
\end{exmp}
\begin{defn}
Let $S_k=\{((x,y),f_{S_k}(x,y)):(x,y)\in X\times Y\}$ be a  two
person soft game and $(x_i,y_j)$, $(x_r,y_s)\in X\times Y$. Then,
Player $k$ is called rational, if the player's soft payoff satisfies
the following conditions:
\begin{enumerate}[a)]
    \item Either $f_{S_k}(x_i,y_j)\supseteq f^k_{X\times
    Y}(x_r,y_s)$ or $f_{S_k}(x_r,y_s)\supseteq f^k_{X\times
    Y}(x_i,y_j)$
    \item If $f_{S_k}(x_i,y_j)\supseteq f^k_{X\times
    Y}(x_r,y_s)$ and $f_{S_k}(x_r,y_s)\supseteq f^k_{X\times
    Y}(x_i,y_j)$, then $f_{S_k}(x_i,y_j)= f^k_{X\times
    Y}(x_r,y_s)$.
\end{enumerate}
\end{defn}

\begin{defn}
Let $S_k=\{((x,y),f_{S_k}(x,y)):(x,y)\in X\times Y\}$ be a  two
person soft game. Then, an action $(x^*,y^*) \in X\times Y$ is
called an optimal action if
$$
 f_{S_k}(x^*,y^*)\supseteq
f_{S_k}(x,y) \textrm{ for all }(x,y) \in X\times Y.
$$
\end{defn}

\begin{defn}
Let $S_k=\{((x,y),f_{S_k}(x,y)):(x,y)\in X\times Y\}$ be a  two
person soft game. Then,
\begin{enumerate}[a)]
\item if $f_{S_k}(x_i,y_j)\supset f_{S_k}(x_r,y_s)$, we says that a player strictly prefers action pair
$(x_i,y_j)$ over action $(x_r,y_s)$,

\item if $f_{S_k}(x_i,y_j)= f_{S_k}(x_r,y_s)$, we says that a player is indifferent between the two
actions,

\item if $f_{S_k}(x_i,y_j)\supseteq f_{S_k}(x_r,y_s)$, we says that a player either prefers $(x_i,y_j)$ to
$(x_r,y_s)$ or is indifferent between the two actions.
\end{enumerate}
\end{defn}

\begin{defn} Let $S_k=\{((x,y),f_{S_k}(x,y)):(x,y)\in X\times Y\}$ be a  two
person soft game for $k=1,2$. Then,
\begin{enumerate}[a)]
    \item If $f_{S_k}(x,y)=\emptyset$ for all $(x,y)\in X\times Y$, then
$S_k$ is called a empty soft set, denoted by $\check{S}_\Phi$.
    \item If $f_{S_k}(x,y)=U$ for all $(x,y)\in X\times Y$, then
$S_k$ is called a full soft game, denoted by $\check{S}_E$.

\end{enumerate}
\end{defn}

Now the two person zero sum game on the classical game theory will
be a two person disjoint game on the soft game theory. It is given
in following definition.

\begin{defn}
A $tps$-game is called a two person disjoint soft game if
intersection of the soft payoff of players is empty set for each
action pairs.

For instance, Example \ref{1} is a two person disjoint soft game.
\end{defn}

\begin{prop}
Let $S_k=\{((x,y),f_{S_k}(x,y)):(x,y)\in X\times Y\}$ be a  two
person disjoint soft game for $k=1,2$. Then,
\begin{description}
\item i.  $({{S_1}^c)}^c= S_1$,
\item ii.  ${{(S_2}^c)}^c= S_2$,
\end{description}
\end{prop}
\begin{prop}
Let $S_k=\{((x,y),f_{S_k}(x,y)):(x,y)\in X\times Y\}$ be a  two
person disjoint soft game for $k=1,2$. Then,
\begin{description}
\item i.  $ S_1 \setminus S_2=S_1$
\item ii.  $ S_2\setminus  S_1=S_2$
\item iii.  $ S_1\cap  S_2=\check{S}_\phi$
\end{description}
\end{prop}
\begin{defn}
A $tps$-game is called a two person universal soft game if union of
the soft payoff of players is universal set for each action pairs.
\end{defn}
For instance, Example \ref{1} is a two person universal soft game.

\begin{prop}
Let $S_k=\{((x,y),f_{S_k}(x,y)):(x,y)\in X\times Y\}$ be a  two
person universal soft game for $k=1,2$. Then,
\begin{description}
\item i.  ${{(S_k}^c)}^c= S_k$, $k=1,2$
\item ii.  $ S_1\cup  S_2=\check{S}_E$
\end{description}
\end{prop}
\begin{prop}
Let $S_k=\{((x,y),f_{S_k}(x,y)):(x,y)\in X\times Y\}$ be a  two
person both universal and disjoint soft game for $k=1,2$. Then,
\begin{description}
\item i.  $ \check{S}_E\setminus S_2=S_1$
\item ii.  $ \check{S}_E\setminus  S_1=S_2$
\end{description}
\end{prop}
\begin{prop}
Let $S_k=\{((x,y),f_{S_k}(x,y)):(x,y)\in X\times Y\}$ be a  two
person both universal and disjoint soft game for $k=1,2$. Then,
\begin{description}
\item i.  ${({S_1}^c)}^c= S_1$,
\item ii.  ${{(S_2}^c)}^c= S_2$,
\item i.  ${{S_1}^c}= S_2$,
\item ii.  ${{S_2}^c} = S_1$,
\end{description}
\end{prop}
\begin{prop}
Let $S_k=\{((x,y),f_{S_k}(x,y)):(x,y)\in X\times Y\}$ be a  two
person both universal and disjoint soft game for $k=1,2$. Then,
\begin{description}
\item i.  $ S_1\setminus S_2=S_1$
\item ii.  $ S_2\setminus  S_1=S_2$
\item iii.  $ S_1\cap  S_2=\check{S}_\phi$
\item iv.  $ S_1\cup  S_2=\check{S}_E$
\end{description}
\end{prop}
\begin{defn}
Let $f_{S_k}$ be a soft payoff function of a $tps$-game ${S_k}$. If
the following properties hold
\begin{enumerate}[a)]
\item $\bigcup_{i=1}^m f_{S_k}(x_i,y_j)=f_{S_k}(x,y)$
\item $\bigcap_{j=1}^n f_{S_k}(x_i,y_j)=f_{S_k}(x,y)$
\end{enumerate}
then $f_{S_k}(x,y)$ is called a soft saddle point of Player k's in
the $tps$-game.
\end{defn}
Note that if $f_{S_1}(x,y)$ is a soft saddle point of a $tps$-game
${S_1}$, then Player 1 can then win at least by choosing the
strategy $x\in X$, and Player 2 can keep her/his loss to at most
$f_{S_1}(x,y)$ by choosing the strategy $y\in Y$. Hence the soft
saddle poind is a value of the $tps$-game.
\begin{exmp}\label{10}
Let $U= \{u_1,u_2,u_3, u_4,u_5,u_6, u_7, u_8,u_9,u_{10}\}$ be  a set
of alternatives,  $X=\{x_1, x_2, x_3,x_4\}$ and $Y=\{y_1, y_2,
y_3\}$ be the strategies for Player 1 and 2, respectively. Then,
$tps$-game of Player 1 is given as in Table 4,
$$
\begin{tabular}{|c|c|c|c|c|c|c|c|}
\hline
${S_1}$& \(y_{1}\) & \(y_{2}\) &  \(y_{3}\)\\
\hline
\(x_1\) & \(\{u_2, u_4, u_7 \} \) & \(\{u_4\}\) &  \( \{ u_4\}\) \\
\hline
\(x_2\) & \(\{u_5\} \) & \(\{u_7\}\) &  \(\{u_4, u_7\} \) \\
\hline
\(x_3\) & \(\{u_2, u_4,u_5, u_7, u_8,u_{10}\} \) & \(\{u_4,u_8\}\) &  \(\{ u_7, u_8\} \) \\
\hline
\(x_4\) & \(\{u_2, u_4,u_5, u_7, u_8\} \) & \(\{u_1,u_4, u_7,u_8\}\) &  \(\{u_4, u_7, u_8\} \) \\
\hline
\end{tabular}
$$
\begin{center}
\footnotesize{\emph{Table 4}}
\end{center}
Clearly,
$$
\begin{array}{l}
\bigcup_{i=1}^4 f_{S_1}(x_i,y_1)= \{u_2, u_4,u_5, u_7, u_8,u_{10}\},
\\
\bigcup_{i=1}^4 f_{S_1}(x_i,y_2)= \{u_1,u_4, u_7,u_8 \},
\\
 \bigcup_{i=1}^4 f_{S_1}(x_i,y_3)=\{u_4,u_7,u_8\},
\end{array}
$$
and
$$
\begin{array}{l}
\bigcap_{j=1}^3 f_{S_1}(x_1,y_j)=\{u_4\},
\\
\bigcap_{j=1}^3 f_{S_1}(x_2,y_j)=\phi,
\\
\bigcap_{j=1}^3 f_{S_1}=\{u_8\},
\\
\bigcap_{j=1}^3 f_{S_1}(x_4,y_j)=\{u_4,u_7,u_8\}.
\end{array}
$$
Therefore, $\{u_4, u_7, u_8\}$ is a soft saddle point of the
$tps$-game, since the intersection of the forth row is equal to the
union of the third column. So, the value of the $tps$-game is
$\{u_4, u_7, u_8\}$.
\end{exmp}
Note that every $tps$-game has not a soft saddle point. (For
instance, in the above example, if $\{u_4,u_7,u_8\}$ is replaced
with $\{u_4,u_7,u_8,u_9\}$ in soft payoff $f_{S_1}(x_4,y_3)$, then a
soft saddle point of the game can not be found.) Saddle point can
not be used for a $tps$-game, soft upper and soft lower values of
the $tps$-game may be used is given in the following definition.
\begin{defn}
Let $f_{S_k}$ be a soft payoff function of a $tps$-game ${S_k}$.
Then,
\begin{enumerate}[i.]
\item
Soft upper value of the $tps$-game, denoted $\overline{v}$, is
defined by
$$
\overline{v}={\cap}_{y\in Y}({\cup}_{x\in X}(f_{S_k}(x,y)))
$$
\item
Soft lower value of the $tps$-game, denoted $\underline{v}$, is
defined by
$$
\underline{v}={\cup}_{x\in X}({\cap}_{y\in Y}(f_{S_k}(x,y)))
$$
\item
If soft upper and soft lower value of a $tps$-game are equal, they
are called value of the $tps$-game, noted by $v$. That is
$v=\underline{v}=\overline{v}$.
\end{enumerate}
\end{defn}
\begin{exmp}\label{199}Let us consider Table 4 in Example \ref{10}.
It is clear that soft upper value $\overline{v}=\{u_4, u_7, u_8\}$
and soft lower value $\underline{v}=\{u_4, u_7, u_8\}$, hence $
\underline{v}=\overline{v}$. It means that value of the $tps$-game
is $\{u_4, u_7, u_8\}$.
\end{exmp}
\begin{ttheorem} $\underline{v}$ and $\overline{v}$ be a soft lower and
soft upper value of a $tps$-game, respectively. Then, the soft lower
value is subset or equal to the soft upper value, that is,
$$\underline{v} \subseteq \overline{v}$$
\end{ttheorem}
Proof:  Assume that $\underline{v}$ be a soft lower value,
$\overline{v}$ be a soft upper value of a $tps$-game and
$X=\{x_1,x_2,...,x_m\}$ and $Y=\{y_1,y_2,...,y_n\}$ are sets of the
strategies for Player 1 and 2, respectively.

We choose $x_i^* \in X$ and $y_j^* \in Y$. Then,
$$
\begin{array}{rl}
\underline{v}&= {\cup}_{x\in X}({\cap}_{y\in Y}(f_{X\times
Y}(x,y)))\\& \subseteq {\cap}_{y\in Y}(f_{X\times Y}(x^*,y))\\&
\subseteq f_{X\times Y}(x^*,y^*)\\&
\subseteq {\cup}_{x\in
X}(f_{X\times Y}(x,y^*))\\& \subseteq
{\cap}_{y\in Y}({\cup}_{x\in X}(f_{X\times Y}(x,y))) \\
\end{array}
$$
i.e.:
$$
\underline{v}={\cup}_{x\in X}({\cap}_{y\in Y}(f_{X\times
Y}(x,y)))\subseteq \overline{v}={\cap}_{y\in Y}({\cup}_{x\in
X}(f_{X\times Y}(x,y)))
$$
proof is valid.
\begin{exmp}\label{190}Let us consider soft upper value $\overline{v}$ and soft lower value $\underline{v}$ in Example \ref{199}. It is clear
that $\overline{v}=\{u_4, u_7, u_8\}\subseteq \underline{v}=\{u_4,
u_7, u_8\}$, hence $ \underline{v}\subseteq \overline{v}$.
\end{exmp}
\begin{ttheorem}
Let $f_{S_k}(x,y)$ be a soft saddle point, $\underline{v}$ be a soft
lower value and $\overline{v}$ be a soft upper value of a
$tps$-game. Then,
$$\underline{v}\subseteq f_{S_k}(x^*,y^*)\subseteq \overline{v}$$
\end{ttheorem}
\textbf{Proof:}  Assume that $f_{S_k}(x^*,y^*)$ be a soft saddle
point, $\underline{v}$ be a soft lower value, $\overline{v}$ be a
soft upper value of a $tps$-game and $X=\{x_1,x_2,...,x_m\}$ and
$Y=\{y_1,y_2,...,y_n\}$ are sets of the strategies for Player 1 and
2, respectively.

We choose $x_i^* \in X$ and $y_j^* \in Y$. Then,

Since $f_{S_k}(x^*,y^*)$ is a soft saddle point, we have
$$\bigcup_{i=1}^m f_{S_k}(x_i,y_j)=\bigcap_{j=1}^n f_{S_k}(x_i,y_j)=f_{S_k}(x^*,y^*)$$

Clearly,
\begin{equation}\label{soft point 1}
eV_L={\cup}_{x\in X}({\cap}_{y\in Y}(f_{X\times Y}(x,y))) \subseteq
\bigcup_{i=1}^m f_{S_k}(x_i,y_j)=f_{S_k}(x^*,y^*)
\end{equation}
 and
\begin{equation}\label{soft point 2}
f_{S_k}(x^*,y^*)=\bigcap_{j=1}^n f_{S_k}(x_i,y_j) \subseteq
eV_U={\cap}_{y\in Y}({\cup}_{x\in X}(f_{X\times Y}(x,y))
\end{equation}

Then, from (\ref{soft point 1}) and (\ref{soft point 2})
$$eV_L\subseteq f_{X\times Y}(x,y)\subseteq eV_U$$
proof is valid.
\begin{ccorollary}
Let $f_{S_k}(x,y)$ be a soft saddle point, $\underline{v}$ be a soft
lower value and $\overline{v}$ be a soft upper value of a
$tps$-game. If $v=\underline{v}=\overline{v}$, then $f_{S_k}(x,y)$
is exactly $v$.
\end{ccorollary}
\begin{exmp}\label{19}Let us consider Table 4 in Example \ref{10} and soft upper
 value $\overline{v}$ and soft lower value $\underline{v}$ in Example \ref{199}.
It is clear that soft saddle point $f_{S_k}(x,y)$ is exactly
$v=\underline{v}=\overline{v}=\{u_4, u_7, u_8\}$.
\end{exmp}
Note that in every $tps$-game, the soft lower value $\underline{v}$
can not be equals to the soft upper value $\overline{v}$. (For
instance, in the above example, if $\{u_4,u_7,u_8\}$ is replaced
with $\{u_4,u_7,u_8,u_9\}$ in soft payoff $f_{S_1}(x_4,y_3)$, then
the soft lower value $\underline{v}$ can not be equals to the soft
upper value $\overline{v}$.) If in a $tps$-game $\underline{v}\neq
\overline{v}$, then to get the solution of the game soft dominated
strategy may be used. We define soft dominated strategy for
$tps$-game as follows.
\begin{defn}
Let ${S_1}$ be a $tps$-game with its soft payoff function $f_{S_1}$.
Then,
\begin{enumerate}[a)]
\item
a strategy $x_i\in X$ is called a soft dominated to another strategy
$x_r\in X$, if $f_{S_1}(x_i,y)\supseteq f_{S_1}(x_r,y)$ for all $y
\in Y$,
\item
a strategy $y_j\in Y$ is called a soft dominated to another strategy
$y_s\in Y$, if $f_{S_1}(x,y_j)\subseteq f_{S_1}(x,y_s)$ for all $x
\in X$
\end{enumerate}
\end{defn}
By using soft dominated strategy, $tps$-games may be reduced by
deleting rows and columns that are obviously bad for the player who
uses them. This process of eliminating soft dominated strategies
sometimes leads us to a solution of a $tps$-game. Such a method of
solving $tps$-game is called a soft elimination method.

The following $tps$-game can be solved by using the method.
\begin{exmp}Let $U= \{u_1,u_2,u_3, u_4,u_5,u_6, u_7, u_8,u_9,u_{10}\}$ be  a set
of alternatives,  $X=\{x_1, x_2, x_3\}$ and $Y=\{y_1, y_2, y_3\}$ be
the strategies for Player 1 and 2, respectively. Then, $tps$-game of
Player 1 is given as in Table 5,
$$
\begin{tabular}{|c|c|c|c|c|c|c|c|}
\hline
${S_1}$& \(y_{1}\) & \(y_{2}\) &  \(y_{3}\)\\
\hline
\(x_1\) & \(\{u_2, u_4, u_7 \} \) & \(\{u_4\}\) &  \( \{ u_4\}\) \\
\hline
\(x_2\) & \(\{u_5\} \) & \(\{u_7\}\) &  \(\{u_4, u_7\} \) \\
\hline
\(x_3\) & \(\{u_2, u_4,u_5, u_7, u_8,u_{10}\} \) & \(\{u_4, u_7,u_8\}\) &  \(\{u_4, u_7, u_8\} \) \\
\hline
\end{tabular}
$$
\begin{center}
\footnotesize{\emph{Table 5}}
\end{center}
The last column is dominated by the middle column. Deleting the last
column we can obtain Table 6 as:
$$
\begin{tabular}{|c|c|c|c|c|c|c|}
\hline
${S_1}$& \(y_{1}\)  &  \(y_{2}\)\\
\hline
\(x_1\) & \(\{u_2, u_4, u_7 \} \) &   \( \{ u_4\}\) \\
\hline
\(x_2\) & \(\{u_5\} \) &  \(\{u_7\} \) \\
\hline
\(x_3\) & \(\{u_2, u_4,u_5, u_7, u_8,u_{10}\} \) & \(\{u_4, u_7, u_8\} \) \\
\hline
\end{tabular}
$$
\begin{center}
\footnotesize{\emph{Table 6}}
\end{center}
Now, in Table 6, the top row is dominated by the bottom row. (Note
that this is not the case in Table 5). Deleting the top row we
obtain Table 7 as:
$$
\begin{tabular}{|c|c|c|c|c|c|c|}
\hline
${S_1}$& \(y_{1}\)  &  \(y_{2}\)\\
 \hline
\(x_2\) & \(\{u_5\} \) &  \(\{u_7\} \) \\
\hline
\(x_3\) & \(\{u_2, u_4,u_5, u_7, u_8,u_{10}\} \) & \(\{ u_4, u_7, u_8\} \) \\
\hline
\end{tabular}
$$
\begin{center}
\footnotesize{\emph{Table 7}}
\end{center}
In Table 7, Player 1 has a soft dominant strategy $x_3$ so that
$x_2$ is now eliminated. Player 2 can now choose between $y_1$ and
$y_2$ and she/he will clearly choose $y_2$. The solution using the
method is $(x_3,y_2)$, that is, value of the $tps$-game is $\{u_4,
u_7, u_8\}$.
\end{exmp}
Note that the soft elimination method cannot be used for some
$tps$-games which do not have a soft dominated strategies. In this
case, we can use soft Nash equilibrium that is defined as follows.
\begin{defn}
Let ${S_k}$ be a $tps$-game with its soft payoff function $f_{S_k}$
for $k=1,2$. If the following properties hold
\begin{enumerate}[a)]
    \item $f_{S_1}(x^*,y^*)\supseteq f_{S_1}(x,y^*)$ for each $x \in X$
    \item $f_{S_2}(x^*,y^*)\supseteq f_{S_2}(x^*,y)$ for each $y \in Y$
\end{enumerate}
then, $(x^*,y^*) \in X\times Y$ is called a soft Nash equilibrium of
a $tps$-game.
\end{defn}
Note that if $(x^*,y^*) \in X\times Y$ is a soft Nash equilibrium of
a $tps$-game, then Player 1 can then win at least $f_{S_1}(x^*,y^*)$
by choosing strategy $x^* \in X$, and Player 2 can win at least
$f_{S_2}(x^*,y^*)$ by choosing strategy $y^* \in Y$. Hence the soft
Nash equilibrium is an optimal action for $tps$-game, therefore,
$f_{S_k}(x^*,y^*)$  is the solution of the $tps$-game for Player
$k$, $k=1,2$.

Following game, given in Example \ref{7}, can be solved by soft Nash
equilibrium, but it is very difficult to solve by using the others
methods.
\begin{exmp}\label{7}Let $U=
\{u_1,u_2,u_3, u_4,u_5,u_6, u_7, u_8,u_9,u_{10}\}$ be  a set of
alternatives, $X=\{x_1, x_2, x_3\}$ and $Y=\{y_1, y_2, y_3\}$ be the
strategies Player 1 and and 2, respectively. Then, $tps$-game of
Player 1 is given as in Table 8,
$$
\begin{tabular}{|c|c|c|c|c|c|c|c|}
\hline ${S_1}$& \(y_{1}\) & \(y_{2}\) &  \(y_{3}\)\\
\hline
\(x_1\) & \(\{u_1,u_2,u_4, u_7 ,u_8,u_9\} \) & \(\{u_1,u_2,u_4, u_7 ,u_8\}\) &  \( \{u_1,u_2,u_3,u_4, u_7 ,u_8\}\) \\
\hline
\(x_2\) & \(\{u_1,u_2,u_3,u_5\} \) & \(\{u_1,u_4, u_7 ,u_8\}\) &  \(\{u_1,u_2,u_3,u_4,u_5, u_7\} \) \\
\hline
\(x_3\) & \(\{u_2, u_5, u_7, u_8,u_{10}\} \) & \(\{u_2,u_4, u_7 ,u_8\}\) &  \(\{u_4,u_5,u_7, u_8,u_{10} \} \) \\
\hline
\end{tabular}
$$
\begin{center}
\footnotesize{\emph{Table 8}}
\end{center}
and $tps$-game of Player 2 is given as in Table 9,
\end{exmp}
$$
\begin{tabular}{|c|c|c|c|c|c|c|c|}
\hline
${S_2}$& \(y_{1}\) & \(y_{2}\) &  \(y_{3}\)\\
\hline
\(x_1\) & \(\{u_3,u_5,u_6,u_{10} \}\) & \(\{u_3, u_5,u_6,u_9,u_{10} \}\) &  \( \{u_5,u_6,u_9,u_{10}\}\) \\
\hline
\(x_2\) & \(\{u_4,u_6, u_7, u_8,u_9,u_{10} \} \) & \(\{u_2,u_3, u_5,u_6,u_9,u_{10} \}\) &  \(\{u_6,  u_8,u_9,u_{10}\} \) \\
\hline
\(x_3\) & \(\{u_1,u_3, u_4,u_6,u_9 \} \) & \(\{u_1,u_3,u_5,u_6,u_9,u_{10}\}\) &  \(\{u_1,u_2,u_3,u_6,u_9\} \) \\
\hline
\end{tabular}
$$
\begin{center}
\footnotesize{\emph{Table 9}}
\end{center}
From the tables, we have
\begin{enumerate}[a)]
    \item $f_{S_1}(x_1,y_2)\supseteq f_{S_1}(x,y_2)$ for each $x \in X$,
    and
    \item $f_{S_2}(x_1,y_2)\supseteq f_{S_2}(x_1,y)$ for each $y \in Y$
\end{enumerate}
then, $(x_1,y_2) \in X\times Y$ is a soft Nash equilibrium.
Therefore, $f_{S_1}(x_1,y_2)=\{u_1,u_2,u_4, u_7 ,u_8\}$ and
$f_{S_2}(x_1,y_2)=\{u_3, u_5,u_6,u_9,u_{10} \}$ are the solution of
the $tps$-game for Player $1$ and Player 2, respectively.
\section{An Application}
In this section, we give a financial problem that are solved by
using both soft dominated strategy and soft saddle point methods.

There are two companies, say Player 1 and Player 2, who
competitively want to increase sale of produces in the country.
Therefore, they give advertisements. Assume that two companies have
a set of different products $U=\{u_1, u_2, u_3, u_4,
u_5,u_6,u_7,u_8\}$ where for $i=1,2,...,8$, the product $u_i$ stand
for ``oil", ``salt",, ``honey", `` jam", `` cheese", ``sugar",
``cooker", and ``jar", respectively. The products can be
characterized by a set of strategy  $X=Y=\{x_i:i=1,2,3\}$ which
contains styles of advertisement where for $j=1,2,3$, the strategies
$x_j$ stand for ``TV", ``radio"  and `` newspaper", respectively.

Suppose that $X=\{x_1, x_2, x_3\}$ and $Y=\{y_1=x_1, y_2=x_2,
y_3=x_3\}$ are strategies of Player 1 and 2, respectively. Then, a
$tps$-game of Player 1 is given as in Table 10.
$$
\begin{tabular}{|c|c|c|c|c|c|c|c|}
\hline
${S_1}$& \(y_{1}\) & \(y_{2}\) &  \(y_{3}\)\\
\hline
\(x_1\) & \(\{u_1,u_2,u_3,u_5,u_8\} \) & \(\{u_1,u_2,u_3,u_4,u_5,u_8\} \) &  \(\{u_3\} \) \\
\hline
\(x_2\) & \(\{u_1,u_3,u_7\} \) & \(\{u_1,u_2,u_3,u_5,u_6,u_7\} \) &  \(\{u_2,u_3\}\) \\
\hline
\(x_3\) & \(\{u_1,u_2,u_3,u_4,u_5\} \) & \(\{u_1,u_2,u_3,u_4,u_5,u_6,u_8\} \) &  \(\{u_1,u_2,u_3\}\) \\
\hline
\end{tabular}
$$
\begin{center}
\footnotesize{\emph{Table 10}}
\end{center}

In Table 10, let us explain action pair $(x_1,y_1)$; if Player 1
select $x_1="TV"$ and Player 2 select $y_1="TV"$, then the soft
payoff of Player 1 is a set $\{u_1,u_2,u_3,u_5,u_8\}$, that is,
$$
f_{S_1}(x_1,y_1)=\{u_1,u_2,u_3,u_5,u_8\}
$$

In this case, Player 1 increase sale of $\{u_1,u_2,u_3,u_5,u_8\}$
and Player 2 decrease sale of $\{u_1,u_2,u_3,u_5,u_8\}$.

We can now solve the game. It is seen in Table 10,
$$
\begin{array}{rcl}
\{u_1,u_2,u_3,u_5,u_8\}         & \subseteq & \{u_1,u_2,u_3,u_4,u_5,u_8\} \\
\{u_1,u_3,u_7\}             & \subseteq & \{u_1,u_2,u_3,u_5,u_6,u_7\}\\
\{u_1,u_2,u_3,u_4,u_5\} & \subseteq &\{u_1,u_2,u_3,u_4,u_5,u_6,u_8\}
\end{array}
$$
the middle column is dominated by the light column. We then deleting
the middle column we obtain Table 11.
$$
\begin{tabular}{|c|c|c|c|c|c|c|c|}
\hline
${S_1}$&  \(y_{1}\) &  \(y_{3}\)\\
\hline
\(x_1\) & \(\{u_1,u_2,u_3,u_5,u_8\}\) &  \(\{u_3\} \) \\
\hline
\(x_2\)  & \(\{u_1,u_3,u_7\}\) &  \(\{u_2,u_3\}\) \\
\hline
\(x_3\)  & \(\{u_1,u_2,u_3,u_4,u_5\}\) &  \(\{u_1,u_2,u_3\} \) \\
\hline
\end{tabular}
$$
\begin{center}
\footnotesize{\emph{Table 11}}
\end{center}
In Table 11, there is no another soft dominated strategy, we can use
soft saddle point method.
$$
\begin{array}{rcl}
\bigcup_{i=1}^3 f_{S_1}(x_i,y_1) & = &\{u_1, u_2, u_3,
u_4,u_5,u_7,u_8\}
\\~~\\
\bigcup_{i=1}^3 f_{S_1}(x_i,y_3)& = & \{u_1, u_2,u_3\}
\\~~\\
\bigcap_{j=1,3} f_{S_1}(x_1,y_j)& = &\{u_3\}
\\~~\\
\bigcap_{j=1,3} f_{S_1}(x_2,y_j)& = &\{u_3\}
\\~~\\
\bigcap_{j=1,3} f_{S_1}(x_3,y_j)& = & \{u_1, u_2,u_3\}
\end{array}
$$
Here, optimal strategy of the game is $(x_3,y_3)$ since
$$
\bigcup_{i=1}^3 f_{S_1}(x_i,y_3)=\bigcap_{j=1,3} f_{S_1}(x_3,y_j)
$$
Therefore, value of the $tps$-game is $ \{u_1, u_2, u_3\}$.

\section{$n$-Person Soft Games}
In many applications the soft games can be often played between more
than two players. Therefore, $tps$-games can be extended to
$n$-person soft games.
\begin{defn}Let $U$ be  a set of alternatives,  $P(U)$ be the power set of $U$
and $X_k$ is the set of strategies of Player $k$, $(k=1,2,...,n)$.
Then, for each Player $k$, an $n$-person soft game ($nps$-game) is
defined by a soft set over $U$ as
$$
S^n_k=\{((x_1,x_2,...,x_n),f_{S^n_k}(x_1,x_2,...,x_n)):(x_1,x_2,...,x_n)\in
X_1\times X_2\times ...\times X_n\}
$$
where $f_{S^n_k}$ is a soft payoff function of Player $k$.

The $nps$-game is played as follows: at a certain time Player 1
chooses a strategy $x_1\in X_1$ and simultaneously each Player $k$
$(k=2,...,n)$ chooses a strategy $x_k\in X_k$ and once this is done
each player $k$ receives the soft payoff
$f_{S^n_k}(x_1,x_2,...,x_n)$.
\end{defn}
\begin{defn}Let $S^n_k=\{((x_1,x_2,...,x_n),f_{S^n_k}(x_1,x_2,...,x_n)):(x_1,x_2,...,x_n)\in
X_1\times X_2\times ...\times X_n\}$ be an $nps$-game. Then, a
strategy $x_k\in X_k$ is called a soft dominated to another strategy
$x\in X_k$, if
$$
f_{S^n_k}(x_1,...,x_{k-1},x_k,x_{k+1},...,x_n) \supseteq
f_{S^n_k}(x_1,...,x_{k-1},x,x_{k+1},...,x_n)
$$
for each strategy $x_i\in X_i$ of player $i$
$(i=1,2,...k-1,k+1,...,n)$, respectively.
\end{defn}
\begin{defn}Let $S^n_k=\{((x_1,x_2,...,x_n),f_{S^n_k}(x_1,x_2,...,x_n)):(x_1,x_2,...,x_n)\in
X_1\times X_2\times ...\times X_n\}$ be an $nps$-game. If for each
player $k$ (k=1,2,...,n) the following properties hold
$$
\begin{array}{rl}
f_{S^n_k}(x^*_1,...,x^*_{k-1},x^*_k,x^*_{k+1},...,x^*_n)\supseteq f_{S^n_k}(x^*_1,...,x^*_{k-1},x,x^*_{k+1},...,x^*_n) \\
\end{array}
$$ for each $x \in X_k$, then $(x^*_1,x^*_2,...,x^*_n) \in S^n_k$ is called a soft Nash
equilibrium of an $nps$-game.
\end{defn}

\section{Conclusion}
In this paper, we first present the basic definitions and results of
soft set theory. We then construct $tps$-games with soft payoffs. We
also give four solution methods for the $tps$-games with examples.
To applied the game to the real world problem we give an example
which shows the methods can be successfully applied to a financial
problem. Finally, we extended the two person soft games to n-person
soft games. The soft games may be applied to many fields and more
comprehensive in the future to solve the related problems, such as;
computer science, decision making, and so on.


\end{document}